\documentclass[]{article}
\newcommand{\pv}{\noindent{\bf Preuve : \/}}
\def\qed{\quad\raise -2pt\hbox{\vrule\vbox
to 10pt{\hrule width 4pt \vfill\hrule}\vrule}\medskip}
\newenvironment{ex}{\smallskip\noindent{\bf Exemple.
}}{\smallskip}
\newenvironment{rmq}{\smallskip\noindent{\bf Remarque.
}}{\smallskip}

\newtheorem{defi}{D\'efinition}

\newcommand\rob{\mathop{\rm rob}}
\newcommand\rb{\mathop{\rm rb}}
\newcommand\ros{\mathop{\rm ros}}
\newcommand\rcb{\mathop{\rm rcb}}
\newcommand\lcb{\mathop{\rm lcb}}
\newcommand\lob{\mathop{\rm lob}}
\newcommand\rcs{\mathop{\rm rcs}}
\newcommand\lcs{\mathop{\rm lcs}}

\newcommand\los{\mathop{\rm los}}

\newcommand\mak{\mathop{\rm mak}}
\newcommand\lmak{\mathop{\rm lmak}}
\newcommand\bmaj{\mathop{ \rm bmaj}}
\newcommand\binv{\mathop{ \rm binv}}

\newcommand\linv{\mathop{ \rm leb}}
\newcommand\rinv{\mathop{ \rm res}}
\newcommand\nrinv{\mathop{ \rm reb}}

\newcommand\stat{\mathop{ \rm stat}}
\newcommand\lb{\mathop{ \rm lb}}
\newcommand\rs{\mathop{ \rm rs}}
\newcommand\ls{\mathop{ \rm ls}}
\newcommand\I{\mathop{\rm I}}
\newcommand\inv{\mathop{\rm inv}}

\def\O{{\cal O}}
\def\S{{\cal S}}
\def\F{{\cal F}}

\def\OP{{\cal O}{\cal P}}
\def\P{{\cal P}}

\newtheorem{prop}{Proposition}

\newtheorem{lem}{Lemme}
\newtheorem{theo}{Th\'eor\`eme}

\newtheorem{conj}{Conjecture}
\title{Nouvelles statistiques de partitions pour
 les $q$-nombres de Stirling de seconde esp\`ece}
\author{G\'erald Ksavrelof et Jiang Zeng\\
\small Institut Girard Desargues\\ \small Universit\'e Claude
Bernard (Lyon 1)\\ \small 69622 Villeurbanne Cedex, France\\
\small \texttt{ksavrelo@desargues.univ-lyon1.fr,
zeng@desargues.univ-lyon1.fr}}
\date{}
\begin{document}
\maketitle

\newdimen\totalhsize
\newdimen\totalvsize
\vbadness=10000 \hoffset=-.45truecm \totalhsize=6.5truein
\totalvsize=9.5truein \textheight=8.25truein \textwidth=5.25truein
\topmargin=0.25truein \headsep=16truept \topskip=11truept
\headheight=11truept \evensidemargin=.5truein
\oddsidemargin=.75truein \footskip=2truepc
\parskip=0truept

\newenvironment{resume}{%
         \small
        \begin{center}%
          {\bfseries {R\'esum\'e}\vspace{-.5em}}%
                 \end{center}%
        \quotation}

%\thispagestyle{empty}%empeche la numerotation de la premiere page
%% pas de numerotation de pages:
%\nopagenumbers

%% et les dimensions suivantes pour le texte:

\hsize = 5.25truein \vsize = 8.25truein

\begin{abstract}
Steingr\'{\i}msson\cite{St} has recently introduced a partition
analogue
of Foata-Zeilberger's mak statistic for permutations and  conjectured 
that its generating function is equal to the classical $q$-Stirling
numbers of second kind.
In this paper we prove a generalization of 
Steingr\'{\i}msson's conjecture~\cite[Conj. 12]{St}.
\end{abstract}
\begin{resume}
Steingr\'{\i}msson\cite{St} a r\'ecemment introduit un analogue en
partitions de la statistique $\mak$  de
 Foata-Zeilberger pour les  permutations et conjectur\'e que leur fonction
g\'en\'eratrice est \'egale aux 
$q$-nombres de Stirling
 de  seconde esp\`ece. Dans cet article nous d\'emontrons 
une g\'en\'eralisation de la conjecture de 
 Steingr\'{\i}msson~\cite[Conj. 12]{St}.
\end{resume}
%%%%%%%%%%%%%%%%%%%%%%%%%%%%%%%%%%%%%%%%%%%%%%%%%%%%%%%%%%%%%%%%%%%%%%%%%%%%%%%
\section{Introduction}
Les  $q$-\emph{nombres de Stirling de seconde esp\`ece}, not\'es
$S_{q}(n,k)$, sont d\'efinis 
par la relation de r\'ecurrence~:
\begin{equation}\label{eq:stirling}
S_{q}(n,k)=\left\lbrace\begin{array}{ll}
q^{k-1}S_{q}(n-1,k-1)+[k]_qS_{q}(n-1,k) \quad&\textrm{si}\quad
1\leq k\leq n,\\ \delta_{n\,k}\quad&\textrm{si}\quad
n=0\,\textrm{ou}\,k=0,
       \end{array}
\right.
\end{equation}
o\`u $[k]_q=1+q+\cdots +q^{k-1}$ pour tout entier $k\geq 1$. 
Carlitz~\cite{CA} et Gould~\cite{GO} ont  \'etudi\'e  pour la premi\`ere fois
ces nombres sous la forme ${\tilde S}_q(n,k)=q^{-{k\choose 2}}S_{q}(n,k)$.
Dans les derni\`eres ann\'ees, beaucoups d'auteurs
ont cherch\'e des interpr\'etations des
$q$-nombres de Stirling ${\tilde S}_q(n,k)$ et  ${S}_q(n,k)$ dans
 diff\'erents
mod\`eles tels que les {\sl partitions}, 
{\sl fonctions \`a croissance restreinte}~\cite{Sa,St,Wa,WW,Wh}, 
{\sl placement de tours} \cite{Mi,GR}, 0-1 {\sl tableaux}~\cite{LE,DL} et 
{\sl juggling patterns}~\cite{EH}.

Une \emph{partition} en $k$ blocs de $[n]=\{1,\, 2,\ldots, n\}$ 
sera not\'e $\pi=B_{1}-\cdots-B_{k}$, o\`u 
$B_{1}, \ldots, B_{k}$ sont les blocs
class\'es par \emph{ordre croissant} de leurs plus petits
\'el\'ements. On note $\P_{n}^{k}$ l'ensemble des
\emph{partitions} en $k$ \emph{blocs} de $[n]$. 
Etant donn\'ee une partition $\pi$ de $[n]$, on r\'epartit
les entiers de $[n]$ en quatre  types  de la mani\`ere suivante~:
\begin{itemize}
\item un \emph{ouvrant} est le plus petit \'el\'ement d'un bloc de $\pi$ ;
\item un \emph{fermant} est le plus grand \'el\'ement d'un bloc de $\pi$ ;
\item un \emph{passant} est un \'el\'ement ni ouvrant ni fermant
d'un bloc de $\pi$ non r\'eduit \`a un seul \'el\'ement ;
\item un \emph{singleton} est l'\'el\'ement d'un bloc de $\pi$ qui n'a qu'un
seul \'el\'ement.
\end{itemize}
L'ensemble des ouvrants, fermants, passants et
singletons de $\pi$ sera not\'e  respectivement par $\O(\pi)$, ${\F}(\pi)$, $\P(\pi)$ et
$\S(\pi)$. Il est \'evident que $\S(\pi)=\O(\pi)\cap {\F}(\pi)$. 

Rappelons qu'un mot 
$w\in [k]^n$ est une {\sl fonction \`a croissance restreinte}
s'il satisfait les conditions suivantes~:
$$
w_1=1 \qquad \textrm{et}\qquad w_i\leq \max\{w_j\,:\, 1\leq j<i\}+1\qquad \textrm{pour tout}\quad i\in [n].
$$
A toute  partition $\pi=B_1-\cdots -B_k\in \P_{n}^{k}$ 
on peut associer une {\sl fonction \`a croissance restreinte}
 $w(\pi)=w_1w_2\ldots w_n$ 
o\`u $w_i$ est  l'indice du bloc de $\pi$ contenant l'entier
$i$ pour $i\in [n]$.

\begin{ex} Si $\pi = 1\,4\,8-2-3\,7\,9-5\,6$, alors
${\cal O}(\pi) =\{1,2,3,5\}$, ${\cal{F}}(\pi)=\{8,2,9,6\}$,
$P(\pi)=\{4,7\}$ et $S(\pi)=\{2\}$ et  $w(\pi)=1\:2\:3\:1\:4\:4\:3\:1\:3$.
\end{ex}

Suivant  Steingr\'{\i}msson~\cite{St} on d\'efinit
les huit statistiques coordonn\'ees sur $\P_{n}^{k}$ comme suit~:
\begin{eqnarray*}
\ros{}_{i}(\pi)&=& \#\{j \in {\cal{O}}(\pi)\,|\,i>j,\,
w_j>w_i\}, \\ 
\rob{}_{i}(\pi)&=& \# \{ j \in
{\O}(\pi)\,|\,i<j,\, w_j>w_i\}, \\
 \rcs{}_{i}(\pi)&=&
\# \{j \in {\cal{F}}(\pi)\,|\,i>j,\, w_j>w_i\},
\\ \rcb{}_{i}(\pi)&=& \# \{j \in {\cal{F}}(\pi)\,|\,i<j,
w_j>w_i\},\\
 \los{}_{i}(\pi)&=& \# \{j \in
{\cal{O}}(\pi)\,|\,i>j, w_j<w_i\},\\
\lob{}_{i}(\pi)&=& \#\{j \in {\cal{O}}(\pi)\,|\,i<j,
w_j<w_i\}, \\ 
\lcs{}_{i}(\pi)&=& \# \{j \in
{\cal{F}}(\pi)\,|\,i>j, w_j<w_i\}, \\
\lcb{}_{i}(\pi)&=& \# \{j \in {\cal{F}}(\pi)\,|\,i<j,
w_j<w_i\}.
\end{eqnarray*}
On d\'efinit ensuite  les huit statistiques 
$\ros$, $\rob$, $\rcs$, $\rcb$, $\lob$, $\los$,  $\lcs$ et
$\lcb$ comme la somme de leurs coordonn\'ees, par exemple, $\ros(\pi)=\sum_{i}\ros_{i}(\pi)$. 
\medskip

\begin{rmq}
$\ros$ est l'abr\'eviation en anglais pour ``\emph{right, opener,
smaller}'', de m\^eme $\lcb$ est celle pour ``\emph{left,
closer, bigger}'', etc (voir \cite{St}).
Certaines de ces statistiques ont \'et\'e introduites dans la lit\'erature 
sous diff\'erentes formes, comme  
les statistiques $lb,\, ls,\, rb$ et $rs$ de \cite[section 2]{WW},
les statistiques ${\I}^M$ et ${\I}^m$ de \cite[section 4]{Mi} et la statistique $inv$ de 
\cite[section 4]{Sa}. Plus exactement, on a les relations suivantes~:
\begin{eqnarray*}
\ros&=&\lb={\I}^M=\inv,\qquad \lcb=\rs,\\
\los&=&\ls={\I}^m,\hskip 1.7cm
\rcb=\rb.
\end{eqnarray*}
D'autre part on voit facilement que la statistique 
$\lob(\pi)\equiv 0$ si les blocs
de $\pi$ sont class\'es par ordre croissant de leurs plus petits \'el\'ements.
Elle ne sera utile que si les blocs de $\pi$ sont class\'es dans un ordre
arbitraire.
\end{rmq}

Inspir\'e par la statistique $\mak$ de Foata et Zeilberger~\cite{FZ}
sur les permutations, Steingr\'{\i}msson~\cite{St} a introduit
les analogues en partitions suivants~:
\begin{defi} Pour tout $\pi\in \P_n^k$, on pose
\begin{eqnarray*}
\mak(\pi)&=& \ros (\pi)+\lcs (\pi),\qquad \lmak{}'(\pi)=
n(k-1)-[\lcb(\pi)+\rob(\pi)],\\ {\mak}'(\pi)&=& \lob (\pi)+\rcb
(\pi),\qquad \lmak{}(\pi)= n(k-1)-[\los(\pi)+\rcs(\pi)].
\end{eqnarray*}
\end{defi}
Donnons un exemple de calculs des statistiques pr\'ec\'edemment
d\'efinies.

\begin{ex}
 Soit  $\pi = 1\,4\,8-2\,9-3\,7-5\,6$, alors on a~:
$$
\begin{array}{ll}
\pi= &  1\,\,4\,\,8-2\,\,9-3\,\,7-5\,\,6 \\ \ros_i: &
0\,\,2\,\,3\hspace{0.45cm}0\,\,2 \hspace{0.45cm} 0\,\,1
\hspace{0.45cm} 0\,\,0\\ \lcs_i: & 0\,\,0\,\,0 \hspace{0.45cm}0\,\,1
\hspace{0.45cm} 0\,\,0  \hspace{0.45cm}0\,\,0\\ \lob_i: &
0\,\,0\,\,0 \hspace{0.45cm}0\,\,0  \hspace{0.45cm}0\,\,0
\hspace{0.45cm} 0\,\,0\\ \rcb_i: &
3\,\,3\,\,1\hspace{0.45cm}2\,\,0 \hspace{0.45cm} 1\,\,0
\hspace{0.45cm} 0\,\,0
\end{array}
\qquad
\begin{array}{ll}
\pi=  &1\,\,4\,\,8-2\,\,9- 3\,\,7- 5\,\,6 \\ \lcb_i: &
0\,\,0\,\,0\hspace{0.45cm} 1\,\,0 \hspace{0.45cm} 2\,\,2
\hspace{0.45cm}3\,\,3\\ \rob_i: & 3\,\,1\,\,0  \hspace{0.45cm}2\,\,0
\hspace{0.45cm} 1 \,\,0 \hspace{0.45cm}0\,\,0\\ \los_i: &
0\,\,0\,\,0 \hspace{0.45cm}1\,\,1 \hspace{0.45cm}2\,\,2
\hspace{0.45cm}3\,\,3\\ \rcs_i: & 0\,\,0\,\,2 \hspace{0.45cm} 0\,\,2
\hspace{0.45cm} 0\,\,1 \hspace{0.45cm}0\,\,0
\end{array}
$$ On en d\'eduit  donc d'apr\`es la d\'efinition~: $$
\begin{array}{ll}
\mak(\pi)=8+1=9,& \lmak{}'(\pi)=27-(11+7)=9,\\
\mak{}'(\pi)=0+10=10, &\lmak(\pi)=27-(12+5)=10.
\end{array}
$$
\end{ex}

Les $q$-nombres de Stirling ont \'et\'e obtenus pour la premi\`ere fois comme 
fonction g\'en\'eratrice dans \cite{Mi}, en termes de statistiques $I^M$ et
$I^m$. Puis, {\sl quatre variantes} de ces statistiques ont \'et\'e introduites dans
\cite{WW}. Celles-ci  ont \'et\'e ensuite \'etendues aux {\sl huit variantes} dans \cite{St}, 
vraisemblablement
\'epuisant toutes les possibilit\'es d'exploiter l'indice d'inversion d'une partition. 
Toutes ces extensions 
ont \'et\'e obtenues par modifications triviales de d\'efinitions; 
or, l'\'etude des distributions conjointes, comme dans \cite{WW}, 
ou celle des statistiques m\'elang\'ees,
comme dans \cite{St}, montre que  ces extensions m\'eritent d'\^etre \'etudi\'ees.

On pourrait distinguer deux types de statistiques sur $\P_n^k$ :
celles dont la v\'erification de la
r\'ecurrence~(\ref{eq:stirling}) est facile et celles dont la
v\'erification de la r\'ecurrence~(\ref{eq:stirling}) est
difficile~\cite{WW}. Par exemple, il est facile (voir
\cite{Mi}) de v\'erifier que
 $\sum_{\pi\in \P_{n}^{k}}q^{\los(\pi)}$ satisfait
(\ref{eq:stirling}). D'autre part, il existe des statistiques
ayant pour fonction g\'en\'eratrice $S_{q}(n,k)$, mais la
r\'ecurrence~(\ref{eq:stirling}) s'av\`ere
 plus difficile \`a v\'erifier (voir \cite{WW}).
Dans cet article nous \'etudions quelques nouvelles statistiques du 
dernier type
sur $\P_n^k$ ayant pour fonction g\'en\'eratrice $S_q(n,k)$. En effet,
cet article a \'et\'e motiv\'e par 
la  conjecture suivante de 
Steingr\'{\i}msson~\cite[Conj.~12]{St}~:
\begin{conj}[Steingr\'{\i}msson]
Les quatre statistiques $\mak$, $\lmak$, $\mak{}'$
et ${\lmak}'$ ont  pour fonction g\'en\'eratrice sur $\P_n^k$
les $q$-nombres de Stirling $S_q(n,k)$, c'est-\`a-dire
$$
\sum_{\pi\in \P_n^k}q^{\mak(\pi)}=\sum_{\pi\in
\P_n^k}q^{\mak{}'(\pi)}=\sum_{\pi\in \P_n^k}q^{{\lmak}'(\pi)}
=\sum_{\pi\in \P_n^k}q^{\lmak(\pi)} =S_q(n,k). 
$$
\end{conj}

Il se trouve que
la statistique $\mak{}'$ est  \'egale \`a  la statistique $\rb$ de 
Wachs et White~\cite{WW}, qui avaient \'etablis, 
parmi d'autres,  le r\'esultat suivant~:
\begin{equation}\label{eq:ww}
\sum_{\pi\in
\P_n^k}q^{\mak{}'(\pi)}=S_q(n,k). 
\end{equation}
En s'appuyant sur le r\'esultat (\ref{eq:ww}) de Wachs et White,
on pourrait d\'emontrer 
la conjecture de Steingr\'{\i}msson ci-dessus \`a partir 
des deux  th\'eor\`emes suivants~:
\begin{theo} Il existe une involution $\varphi$ sur $\P_n^k$ telle que
pour tout $\pi\in \P_n^k$,  on a
 $\mak(\pi)=\mak{}'(\varphi(\pi))$.
\end{theo}
\begin{theo} Pour tout  $\pi\in \P_n^k$ on a
$$
\mak(\pi)=\lmak{}'(\pi),\qquad \mak{}'(\pi)=\lmak(\pi).
$$
\end{theo}

En fait, l'approche que  nous proposons dans cet article est 
 ind\'ependante du r\'esultat (\ref{eq:ww}) de Wachs et White et 
 a permis de trouver une nouvelle 
 statistque $\mak_l$ g\'en\'eralisant $\mak$.

\begin{defi}
Soit
 $\pi =B_{1}-\cdots - B_{k}\in \P_{n}^k$ et 
$w(\pi)~=~w_1\ldots w_n$. Pour tout $b\in [n]$ on pose
$\nrinv(b, \pi)=\#\{a\mid w_a>w_b \,\,\textrm{et}\,\, a>b\}$.
Pour tout $ l\in [k]$  
on d\'esigne respectivement
par $p(B_l)$ et $g(B_l)$ le {\sl plus petit} et le {\sl plus grand} \'el\'ement de
$B_l$, et on d\'efinit
$$
{\mak}_{l}(\pi)=\mak(\pi)-\nrinv(g(B_{l}), \pi)+k-l.
$$
\end{defi}
On remarque que
lorsque $l=k$ on retrouve la $\mak$ ordinaire, i.e.  $\mak_{k}=\mak$.

\begin{ex}
Soit  $\pi = 1\,4\,8-2-3\,7\,9-5\,6$, alors on a $g(B_1)=8$,
$g(B_2)=2$, $g(B_3)=9$ et $g(B_4)=6$. Ainsi~:
\begin{center}
\begin{tabular}{ll}
$\nrinv(g(B_1),\pi)=1$ & $\Longrightarrow \qquad \mak_1(\pi)=11$,\\
$\nrinv(g(B_2),\pi)=5$ & $\Longrightarrow \qquad \mak_2(\pi)=6$,\\
$\nrinv(g(B_3),\pi)=0$ & $\Longrightarrow \qquad \mak_3(\pi)=10$,\\
$\nrinv(g(B_4),\pi)=0$ & $\Longrightarrow \qquad \mak_4(\pi)=9$.\\
\end{tabular}
\end{center}
\end{ex}

Le th\'eor\`eme suivant g\'en\'eralise la conjecture de Steingr\'{\i}msson
sur la statistique $\mak$.
\begin{theo} Pour $1\leq l\leq k$, on a
$$ \sum_{\pi \in \P_n^k}q^{{\mak}_{l}(\pi)}=S_q(n,k). $$
\end{theo}

Nous donnons les d\'emonstrations de ces
trois th\'eor\`emes respectivement dans 
les trois sections suivantes et 
terminons l'article avec quelques remarques
sur les probl\`emes ouverts.

%%%%%%%%%%%%%%%%%%%%%%%%%%%%%%%%%%%%%%%%%%%%%%%%%%%%%%%%%%%%%%%%
%%%%%%%%%%%%%%%
\section{Preuve du th\'eor\`eme~1}
%%%%%%%%%%%%%%%%%%%%%%%%%%%%%%%%%%%%%%%%%%%%%%%%%%%%%%%%%%%%%%%%
%%%%%%%%%%%%%%%
Nous avons besoin de quelques d\'efinitions suppl\'ementaires.
Pour tout ensemble fini d'entiers $B$ et entier $i$, on note
$B(\leq i)$ la restriction de $B$ sur $[i]$, qui est soit
\emph{complet}, si l'ensemble $B$ est inclus  dans $[i]$, soit
\emph{incomplet}, si une partie non vide de $B$ est dans $[i]$ et
l'autre partie non vide dans $[n]\setminus [i]$, soit \emph{vide}
si $B\cap [i]=\emptyset$. 

\begin{defi} Soit
$\pi=B_1-B_2-\cdots -B_k$  une partition de $\P_n^k$ et
$T_0=\emptyset$.
Pour $i=1,\ldots, n$,
on d\'efinit  la $i^e$ \emph{trace}  de $\pi$  comme la partition $T_i$ de 
$[i]$~: 
$$ 
T_i=B_1(\leq i)-B_2(\leq i)-\cdots -B_k(\leq i).
$$ 
On note 
le nombre de {\sl blocs incomplets}  dans $T_{i-1}$ par $l_{i}(\pi)$,
 et le nombre  de {\sl blocs incomplets}
situ\'es \`a gauche du bloc contenant $i$ dans  $T_i$ par $\gamma_{i}(\pi)-1$. 
\end{defi}

\begin{ex}
Si $\pi=1\,4\,8-2-3\,7\,9-5\,6$, alors les traces, $l_i(\pi)$
 et $\gamma_i(\pi)$ sont donn\'es par~: 
$$
\begin{array}{lllll}
T_1=1 \cdot &\hspace{1cm}&l_1=0 &\hspace{1cm}&\gamma_1=1\\
T_2= 1 \cdot - 2&\hspace{1cm}&l_2=1&\hspace{1cm}&\gamma_2=2\\
T_3=1 \cdot - 2 - 3 \;\cdot&\hspace{1cm}&l_3=1&\hspace{1cm}&\gamma_3=2\\
T_4=1\;4 \cdot - 2 - 3\; \cdot&\hspace{1cm}&l_4=2&\hspace{1cm}&\gamma_4=1\\
T_5=1 \;4 \cdot - 2 - 3 \cdot - 5 \;\cdot&\hspace{1cm}&l_5=2&\hspace{1cm}&\gamma_5=3\\
T_6=1 \;4 \cdot - 2 - 3 \cdot -5\;6&\hspace{1cm}&l_6=3&\hspace{1cm}&\gamma_6=3\\
T_7=1 \;4 \cdot - 2 - 3 \;7 \cdot - 5\;6&\hspace{1cm}&l_7=2&\hspace{1cm}&\gamma_7=2\\
T_8=1 \;4\;8-2-3 \;7 \cdot - 5\;6&\hspace{1cm}&l_8=2&\hspace{1cm}&\gamma_8=1\\
T_9=1 \;4\;8-2-3 \;7\;9 - 5\;6&\hspace{1cm}&l_9=1&\hspace{1cm}&\gamma_9=1
\end{array}
$$
o\`u on ajoute un {\sl point} dans chaque {\sl bloc incomplet}. 
\end{ex}

Il est clair qu'une partition est enti\`erement
d\'etermin\'ee par ses traces successives  $ T_1, T_2,\ldots, T_n$ ou 
par la suite $(l_1,\, \gamma_1),\ldots, (l_n,\, \gamma_n)$. 
D'autre part,
pour tout $i\in [n]$ on voit que 
\begin{equation}\label{eq:calcul}
k=\left\lbrace\begin{array}{ll} l_i+\#\{a\in {\O}\,|\,
a>i\}+\#\{a\in {\F}\,|\, a<i\},\quad &\textrm{si}\quad {i\in
\P\cup \F_s},\\ 1+l_i+\#\{a\in {\O}\,|\, a>i\}+\#\{a\in {\F}\,|\,
a<i\},\quad &\textrm{si}\quad {i\in \O_s\cup\S}.
\end{array}\right.
\end{equation}
D'o\`u on tire en sommant sur tous les $i$~:
\begin{equation}\label{eq:rel}
nk =\#\O+\sum_{i=1}^n(l_i+\#\{a \in\O\,|\,a>i\}+\#\{a \in\F\,|\,a<i\}). 
\end{equation}
\begin{lem}
Soit $\pi \in \P_n^k$ une partition fix\'ee, on pose
 $\O=\O(\pi)$, $\F_s=\F(\pi)\setminus \S(\pi)$, $\P=\P(\pi)$, $l_i=l_i(\pi)$ et
$\gamma_i=\gamma_i(\pi)$. Alors on a les identit\'es suivantes~:
\begin{eqnarray}
\mak(\pi)&=&\sum_{i\in \F_s\cup \P}(l_i-\gamma_i)
+\sum_{i=1}^n\#\{a\in \F|a<i\},\\
\mak{}'(\pi)&=&\sum_{i\in \F_s\cup \P}(k-\gamma_i)+ 
\sum_{i\in \O}(k-1-l_i)-\sum_{i=1}^n\#\{a\in \F|a<i\},\\
\sum_{i\in {\F_s}}l_i&=&\sum_{i\in {\O_s}}(l_i+1).\label{eq:3}
\end{eqnarray}
\end{lem}
\pv Pout tout $i\in [n]$, comme
$\lcs_{i}(\pi)$ est le nombre de blocs complets \`a gauche du
bloc contenant $i$ dans la $i^e$-trace de $\pi$, on a
$\lcs_{i}(\pi)= \#\{a\in \F(\pi)|a<i\}$;
d'autre part, $\ros_{i}(\pi)$ est le nombre
de blocs (complets ou incomplets) \`a droite du bloc contenant
$i$ dans la $i^e$-trace de $\pi$,  ainsi
$\ros{}_{i}(\pi)=l_{i}-\gamma_i$ si $i \in \F_s\cup \P$ et
$\ros{}_{i}(\pi)=0$ 
si $i \in \O$. 
 D'o\`u  la premi\`ere \'egalit\'e. De m\^eme, 
comme $\rcb_{i}(\pi)$ comptent non seulement les blocs incomplets
\`a droite du bloc contenant $i$ dans $T_i$, c'est \`a dire
$l_{i}-\gamma_i$ si $i\in \F_s\cup \P$ et $0$ sinon, mais aussi les
blocs qui ne sont pas encore cr\'e\'es, soit
$$\left\lbrace\begin{array}{ll}
k-l_{i}-\#\{a\in\F\,|\,a<i\}&\qquad \textrm{si}\quad i\in \F_s\cup \P ,\\
k-l_{i}-\#\{a\in\F\,|\,a<i\}-1&\qquad \textrm{si}\quad i\in \O.
  \end{array}\right.
$$
En sommant sur $i$ on obtient la seconde \'egalit\'e. 
Enfin, on remarque qu'il y a autant de fermants que d'ouvrants
dans $\pi$, et $\forall i \in [n]$, $l_i \geq 0$ et
 $l_{i+1}=l_i+1$ (resp.  $l_{i+1}=l_i-1$~)
si $i\in \O_s$ (resp.  si $i\in \F_s$). On va construire  une
bijection de $\O_s$ dans $\F_s$ telle que si $a\mapsto a'$, alors
$l_{a}+1=l_{a'}$. Ce qui d\'emontre  clairement la troisi\`eme \'egalit\'e. En effet, soit 
$\O_s=\{a_1,\ldots, a_r\}$ tel
que $a_1<a_2<\cdots <a_r$.
 Comme $l_{a_1}=0$  et $l_g=1$, o\`u  $g$
est le plus grand fermant de $\pi$, on d\'efinit $a_1'$ comme le
plus petit fermant tel que  $l_{a_1}+1=l_{a_1'}$. Supposons ainsi
d\'efinis les
 $i-1$ fermants  $a_1',\ldots, a_{i-1}'$ associ\'es avec
les  $i$ premiers ouvrants $a_1,\ldots, a_{i-1}$ respectivement. A
$a_{i}$ on associe le plus petit fermant, soit $a_{i}'$,  dans
$\F_s\setminus\{a_1',\ldots, a_{i-1}'\}$ tel que
 $l_{a_{i}}+1=l_{a_{i}'}$.
\qed

Gr\^ace \`a la notion de trace, on est maintenant en mesure de d\'ecrire
une involution
$\varphi:\P_n^k\,\longrightarrow\,\P_n^k$ d\'efinie par
l'algorithme suivant~:
\begin{enumerate}
  \item Etant donn\'ee une partition $\pi$ de $[n]$, on partage
l'ensemble $[n]$ en quatre  parties $\S(\pi)$,
$\O_s(\pi)=\O(\pi)\setminus\S(\pi)$,
$\F_s(\pi)=\F(\pi)\setminus\S(\pi)$ et $\P(\pi)$, et on calcule
les $\gamma_{i}$ pour tout $i \in [n]$. Notons $f$ (resp. $p$) la
suite croissante des \'el\'ements de $\F_s(\pi)$ (resp.
$\P(\pi)$). On forme alors les deux matrices~: 
$$
\left(\begin{array}{c} f\\ \gamma
\end{array}\right)=\left(\begin{array}{cccc} f_1&f_2&\ldots &f_r\\
\gamma_{f_1}&\gamma_{f_2}&\ldots &\gamma_{f_r}
\end{array}\right),\quad
\left(\begin{array}{c} p\\ \gamma
\end{array}\right)=
\left(\begin{array}{cccc} p_1&p_2&\ldots &p_s\\
\gamma_{p_1}&\gamma_{p_2}&\ldots &\gamma_{p_s}
\end{array}\right).
$$
\item On commence par d\'efinir les quatre ensembles correspondants de $\pi'$~:
\begin{equation}
\begin{array}{lll}
\S'&=\{n+1-i\,|\, i\in \S(\pi)\},\quad  \O'_s&=\{n+1-i\,|\, i\in
\F_s(\pi)\},\\
{\F}'_s&=\{n+1-i\,|\, i\in \O_s(\pi)\},\quad 
\P'&=\{n+1-i\,|\,i\in \P(\pi)\}. 
\end{array}
\end{equation}
On note que 
$$
\O'=\O'_s\cup \S'\quad\textrm{et}\quad \F'=\F_s\cup \S'.
$$
Soient $f'$ et $p'$ les suites
croissantes des \'el\'ements de $\F'_s$ et $\P'$, formons les deux
matrices~: $$ \left(\begin{array}{c} f'\\ \gamma'
\end{array}\right)=\left(\begin{array}{cccc} f_1'&f_2'&\ldots
&f_r'\\ \gamma_{f_1}&\gamma_{f_2}&\ldots &\gamma_{f_r}
\end{array}\right),\quad
\left(\begin{array}{c} p'\\ \gamma'
\end{array}\right)=
\left(\begin{array}{cccc} p_1'&p_2'&\ldots &p_s'\\
\gamma_{p_s}&\gamma_{p_{s-1}}&\ldots &\gamma_{p_1}
\end{array}\right).
$$
\item On commence par construire une partition $\pi_0$ de $\O'$, 
dont les blocs sont tous des singletons. Un singleton $\{i\}$ est dit 
complet (resp. incomplet)  si $i\in \S'$ (resp. sinon).
Soit 
$$ \pi_0=B_1-B_2-\cdots -B_k.$$ 
Supposons ensuite que
$x_1,x_2,\ldots, x_{n-k}$ est le {\sl r\'earrangement croissant} des
\'el\'ements de $\F'_s\cup \P'$. Pour $j=1,\ldots, n-k$ on construit  $\pi_j$ en 
ins\'erant  $x_j$ dans l'un des blocs de $\pi_{j-1}$ de sorte que
$\gamma{}'_{x_j}=\gamma_{x_j}(\pi_j)$. Rappelons qu'un bloc de $\pi_j$ est consid\'er\'e 
complet s'il d\'ebute avec un \'el\'ement de $\O'$ et termine avec un \'el\'ement de $\F'$.
\item D\'efinissons
$\varphi(\pi)=\pi'=\pi_{n-k}$, alors $\S(\pi')=\S'$,
$\O_s(\pi')=\O'_s$, $\F_s(\pi')=\F'_s$, $\P(\pi')=\P'$.
\end{enumerate}
On v\'erifie que $\varphi$ est une involution sur $\P_n^k$ telle
que $\F_s(\varphi(\pi))=\{n+1-i\,|\,i\in \O_s(\pi)\}$.
\medskip

\begin{ex} Prenons la partition $\pi=1\,4\,8-2-3\,7\,9-5\,6$, alors 
$${\O}_s=\{1,3,5\},\quad
{\F}_s~=~\{8,9,6\}, \quad {\P}=\{4,7\}, \quad  \S=\{2\}. 
$$
On en d\'eduit
donc ${\F}_s'=\{9,7,5\}$, ${\O}_s'=\{2,1,4\}$, ${\P}'=\{6,3\}$ et
$\S'=\{8\}$. Ainsi on a $$ \left(\begin{array}{c} f\\ \gamma
\end{array}\right)=\left(\begin{array}{ccc} 6&8&9\\ 3&1&1
\end{array}\right),\quad
\left(\begin{array}{c} p\\ \gamma
\end{array}\right)=
\left(\begin{array}{cc} 4&7\\ 1&2
\end{array}\right),
$$ et puis
$$ \left(\begin{array}{c} f'\\ \gamma'
\end{array}\right)=\left(\begin{array}{ccc} 5&7&9\\ 3&1&1
\end{array}\right),\quad
\left(\begin{array}{c} p'\\ \gamma'
\end{array}\right)=
\left(\begin{array}{cc} 3&6\\ 2&1
\end{array}\right).
$$ On obtient d'abord les $k=4$ blocs (complets ou incomplets)
form\'es des \'el\'ements de $\O'$~: 
$$ 
\pi_0=1\; \cdot - \;2\;\cdot -\; 4\;\cdot -\; 8. 
$$ 
On
ins\`ere  successivement les  \'el\'ements $i$ de $\F_s'\cup \P'$
en tenant compte de $\gamma_i$ et on obtient
\begin{eqnarray*}
\pi_1&=&1\;\cdot -\; 2\; 3 \cdot -\; 4 \cdot - \; 8\\ 
\pi_2&=&1\; \cdot - \;2\; 3 \cdot - \; 4\;5  -\; 8\\ 
\pi_3&=&1\;6\; \cdot - \;2\; 3 \cdot -\; 4\;5 -\; 8\\
 \pi_4&=&1\; 6\; 7 -\; 2\;3 \cdot -\; 4\; 5 - \;8\\ 
\pi_5&=&1\;6\; 7  -\; 2\; 3\; 9 - \;4\;5  - \;8.
\end{eqnarray*}
D'o\`u $\pi'=\pi_5=1\;6\;7\;-\;2\;3\;9\;-\;4\;5\;-\;8$. On v\'erifie que
$(\pi')'=\pi$.
\end{ex}

\begin{rmq} On pourrait d\'ecrire 
l'involution $\varphi$ dans le mod\`ele des {\sl chemins de Motzkin
valu\'es}~\cite{Fl}.
En effet \`a chaque partition $\pi \in \P_n^k$, on peut associer un
chemin $T(\pi)$ d\'efinie comme suit :
\begin{itemize}
\item A chaque  \'el\'ement de $\O_s$ (resp. $\S$), on associe un pas
nord-est (resp. est)
\'etiquet\'e $1$ (resp. $1^*$).
\item  A chaque  \'el\'ement $i \in \F_s$ (resp. $\P$), on associe un pas
sud-est (resp. est)
\'etiquet\'e $\gamma_i$.
\end{itemize}
La correspondance $\pi\longrightarrow \varphi(\pi)=\pi'$ se pr\'esente alors
comme  une symm\'etrie par
rapport \`a l'axe des ordonn\'ees des chemins correspondants
o\`u l'\'etiquettage de $T(\pi')$ est le suivant~:
\begin{itemize}
\item Les pas ``est'' gardent la m\^eme valeur.
\item Les pas ``nord-est'' reste \'etiquet\'es $1$.
\item On \'etiquette les pas ``sud-est'' en reprenant les \'etiquettes
des pas ``sud-est de $T(\pi)$ dans le
m\^eme ordre.
\end{itemize}

\begin{ex}
Soit $\pi=1\,4\,8-2-3\,7\,9-5\,6$ alors :
$$
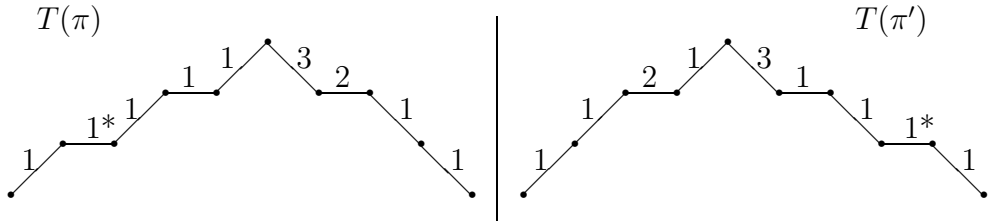
\setlength{\unitlength}{0.17mm}%
\begingroup\makeatletter\ifx\SetFigFont\undefined
% extract first six characters in \fmtname
\def\x#1#2#3#4#5#6#7\relax{\def\x{#1#2#3#4#5#6}}%
\expandafter\x\fmtname xxxxxx\relax \def\y{splain}%
\ifx\x\y   % LaTeX or SliTeX?
\gdef\SetFigFont#1#2#3{%
  \ifnum #1<17\tiny\else \ifnum #1<20\small\else
  \ifnum #1<24\normalsize\else \ifnum #1<29\large\else
  \ifnum #1<34\Large\else \ifnum #1<41\LARGE\else
     \huge\fi\fi\fi\fi\fi\fi
  \csname #3\endcsname}%
\else
\gdef\SetFigFont#1#2#3{\begingroup
  \count@#1\relax \ifnum 25<\count@\count@25\fi
  \def\x{\endgroup\@setsize\SetFigFont{#2pt}}%
  \expandafter\x
    \csname \romannumeral\the\count@ pt\expandafter\endcsname
    \csname @\romannumeral\the\count@ pt\endcsname
  \csname #3\endcsname}%
\fi
\fi\endgroup
\begin{picture}(760,241)(20,400)
\thinlines
\put(400,580){\line( 0,-1){160}}
\put(420,440){\line( 1, 1){ 80}}
\put(500,520){\line( 1, 0){ 40}}
\put(540,520){\line( 1, 1){ 40}}
\put(580,560){\line( 1,-1){ 40}}
\put(620,520){\line( 1, 0){ 40}}
\put(660,520){\line( 1,-1){ 40}}
\put(700,480){\line( 1, 0){ 40}}
\put(740,480){\line( 1,-1){ 40}}
\put( 20,460){\makebox(0,0)[lb]{\smash{\SetFigFont{12}{14.4}{rm}    1   }}}
\put( 70,485){\makebox(0,0)[lb]{\smash{\SetFigFont{12}{14.4}{rm}  1*}}}
\put(100,500){\makebox(0,0)[lb]{\smash{\SetFigFont{12}{14.4}{rm}   1}}}
\put(145,525){\makebox(0,0)[lb]{\smash{\SetFigFont{12}{14.4}{rm}     1}}}
\put(175,540){\makebox(0,0)[lb]{\smash{\SetFigFont{12}{14.4}{rm}     1}}}
\put(235,540){\makebox(0,0)[lb]{\smash{\SetFigFont{12}{14.4}{rm}   3}}}
\put(265,525){\makebox(0,0)[lb]{\smash{\SetFigFont{12}{14.4}{rm}    2}}}
\put(315,500){\makebox(0,0)[lb]{\smash{\SetFigFont{12}{14.4}{rm}       1}}}
\put( 20,440){\line( 1, 1){ 40}}
\put( 60,480){\line( 1, 0){ 40}}
\put(100,480){\line( 1, 1){ 40}}
\put(140,520){\line( 1, 0){ 40}}
\put(180,520){\line( 1, 1){ 40}}
\put(220,560){\line( 1,-1){ 40}}
\put(260,520){\line( 1, 0){ 40}}
\put(300,520){\line( 1,-1){ 80}}
\put(355,460){\makebox(0,0)[lb]{\smash{\SetFigFont{12}{14.4}{rm}     1}}}
\put(680,570){\makebox(0,0)[lb]{\smash{\SetFigFont{12}{14.4}{rm}$T(\pi')$}}}
\put(420,460){\makebox(0,0)[lb]{\smash{\SetFigFont{12}{14.4}{rm}  1}}}
\put(457,500){\makebox(0,0)[lb]{\smash{\SetFigFont{12}{14.4}{rm}      1}}}
\put(505,525){\makebox(0,0)[lb]{\smash{\SetFigFont{12}{14.4}{rm}    2}}}
\put(540,540){\makebox(0,0)[lb]{\smash{\SetFigFont{12}{14.4}{rm}    1}}}
\put(595,540){\makebox(0,0)[lb]{\smash{\SetFigFont{12}{14.4}{rm}     3}}}
\put(625,525){\makebox(0,0)[lb]{\smash{\SetFigFont{12}{14.4}{rm}       1}}}
\put(675,500){\makebox(0,0)[lb]{\smash{\SetFigFont{12}{14.4}{rm}      1}}}
\put(710,485){\makebox(0,0)[lb]{\smash{\SetFigFont{12}{14.4}{rm}  1*}}}
\put(755,460){\makebox(0,0)[lb]{\smash{\SetFigFont{12}{14.4}{rm}     1}}}
\put( 40,570){\makebox(0,0)[lb]{\smash{\SetFigFont{12}{14.4}{rm}$T(\pi)$}}}
\put( 16,438){\makebox(0,0)[lb]{\smash{\SetFigFont{5}{4.4}{rm}$\bullet$}}}
\put( 57,478){\makebox(0,0)[lb]{\smash{\SetFigFont{5}{4.4}{rm}$\bullet$}}}
\put( 97,478){\makebox(0,0)[lb]{\smash{\SetFigFont{5}{4.4}{rm}$\bullet$}}}
\put( 137,518){\makebox(0,0)[lb]{\smash{\SetFigFont{5}{4.4}{rm}$\bullet$}}}
\put( 177,518){\makebox(0,0)[lb]{\smash{\SetFigFont{5}{4.4}{rm}$\bullet$}}}
\put( 217,558){\makebox(0,0)[lb]{\smash{\SetFigFont{5}{4.4}{rm}$\bullet$}}}
\put( 257,518){\makebox(0,0)[lb]{\smash{\SetFigFont{5}{4.4}{rm}$\bullet$}}}
\put( 297,518){\makebox(0,0)[lb]{\smash{\SetFigFont{5}{4.4}{rm}$\bullet$}}}
\put( 337,478){\makebox(0,0)[lb]{\smash{\SetFigFont{5}{4.4}{rm}$\bullet$}}}
\put( 377,438){\makebox(0,0)[lb]{\smash{\SetFigFont{5}{4.4}{rm}$\bullet$}}}
\put( 417,438){\makebox(0,0)[lb]{\smash{\SetFigFont{5}{4.4}{rm}$\bullet$}}}
\put( 457,478){\makebox(0,0)[lb]{\smash{\SetFigFont{5}{4.4}{rm}$\bullet$}}}
\put( 497,518){\makebox(0,0)[lb]{\smash{\SetFigFont{5}{4.4}{rm}$\bullet$}}}
\put( 537,518){\makebox(0,0)[lb]{\smash{\SetFigFont{5}{4.4}{rm}$\bullet$}}}
\put( 577,558){\makebox(0,0)[lb]{\smash{\SetFigFont{5}{4.4}{rm}$\bullet$}}}
\put( 617,518){\makebox(0,0)[lb]{\smash{\SetFigFont{5}{4.4}{rm}$\bullet$}}}
\put( 657,518){\makebox(0,0)[lb]{\smash{\SetFigFont{5}{4.4}{rm}$\bullet$}}}
\put( 697,478){\makebox(0,0)[lb]{\smash{\SetFigFont{5}{4.4}{rm}$\bullet$}}}
\put( 737,478){\makebox(0,0)[lb]{\smash{\SetFigFont{5}{4.4}{rm}$\bullet$}}}
\put( 777,438){\makebox(0,0)[lb]{\smash{\SetFigFont{5}{4.4}{rm}$\bullet$}}}
\end{picture}
$$
A partir du chemin $T(\pi)$
on retrouve alors la partition $\varphi(\pi)=1\,6\,7-2\,3\,9-4\,5-8$ .
\end{ex}
\end{rmq}
Reste \`a montrer que $\mak{}'(\varphi(\pi))=\mak(\pi)$, ce
qui, en vertu du lemme~1,  \'equivaut \`a~:
\begin{eqnarray}
&&\sum_{i\in {\F_s}\cup \P}(l_i-\gamma_i) +\sum_{i=1}^n\#\{a\in
{\F}\,|\,a<i\}\nonumber\\
&&{}= \sum_{i\in {\F'_s}\cup \P'}(k-\gamma_i')+
\sum_{i\in {\O'}}(k-l'_i-1)-\sum_{i=1}^n\#\{a\in {\F'}\,|\,a<i\}.
\end{eqnarray}
Or la construction de $\varphi(\pi)$  exige que les suites $(\gamma_i)_{i\in \F'_s}$ et 
$(\gamma_i)_{i\in \P'_s}$
soient respectivement des  r\'earrangements de $(\gamma_i)_{i\in \F_s}$ et
$(\gamma_i)_{i\in \P_s}$.
Ainsi, compte tenu des d\'efinitions (8) de  $\F'$,
$\O'$, $\P'$, $\S'$, $\F'_s$ et $\O'_s$, l'identit\'e (9)  peut s'\'ecrire
comme suit~:
$$
\sum_{i\in {\F_s}\cup \P}l_i+\sum_{i=1}^n\left(
\#\{a\in {\F}\,|\,a<i\} +\#\{a\in{\O}\,|\,a>i\}\right)
= nk-\sum_{i \in \O'}(l'_i+1).
$$
En appliquant l'\'equation (\ref{eq:rel}), on voit que
l'\'egalit\'e ci-dessus
\'equivaut \`a~:
\begin{equation}\label{eq:ker}
\sum_{i\in {\O}}(l_i+1)=\sum_{i\in {\O'}}(l'_i+1).
\end{equation}
Remarquant que $j\in \O'$ si et seulement si $\bar j=n+1-j\in \F$, on d\'eduit de 
(\ref{eq:calcul}) que  
\begin{eqnarray*}
l_i'&=&k-1-\#\{a\in \O'\mid a>i\}-\#\{a\in \F'\mid a<i\}\\
&=&k-1-\#\{\bar b\in \F\mid \bar b<\bar i\}
-\#\{\bar b\in \O\mid \bar b>\bar i\}.
\end{eqnarray*}
Ce qui  montre  que
$l_i'=l_{\bar i}$ si $i\in \S'$ et $l_i'=l_{\bar i}+1$ 
si $i\in F'_s$.  Il en r\'esulte que l'identit\'e (\ref{eq:ker}) \'equivaut
\`a (\ref{eq:3}). Ce qu'il fallait d\'emontrer.

%%%%%%%%%%%%%%%%%%%%%%%%%%%%%%%%%%%%%%%%%%%%%%%%%%%%%%%
\section{Preuve du th\'eor\`eme~2}
%%%%%%%%%%%%%%%%%%%%%%%%%%%%%%%%%%%%%%%%%%%%%%%%%%%%%%%
On commence par quelques d\'efinitions et notations.
Pour \^etre coh\'erent avec les notations de
Steingr\'{\i}mmson \cite{St},
on utilise r\'espectivement les abr\'eviations
de \emph{ right element smaller}, \emph{right element bigger} et 
\emph{left  element bigger}
 pour les statistiques $\rinv$, $\nrinv$ et $\linv$.
\begin{defi}
Soit $\pi = B_{1}-\cdots - B_{k}$ une partition de $\P_{n}^{k}$ et
$b\in B_{j}$ fix\'e. On d\'efinit d'abord, pour tout  $i>j$, 
$\rinv{}(b,B_{i})=\# \{a \in B_{i} \,|\, b>a\}$,
$\nrinv{}(b,\, B_{i})=\# \{a \in B_{i} \,|\, b<a\}$, 
et pour tout $i<j$, 
$\linv{}(b,B_{i})=\# \{a \in B_{i} \,|\, b<a\}$; et puis
$$
\rinv{}(b,\, \pi)=\sum_{i>j}\rinv{}(b,B_{i}),
\quad 
\nrinv{}(b,\, \pi)=\sum_{i>j}\nrinv{}(b,B_{i}),\quad 
\linv{}(b,\,\pi)=\sum_{i<j}\linv{}(b,B_{i}).
$$
Enfin on note $b_j$ le cardinal de $B_j$  pour tout 
$j\in [k]$ et pose
$$
\linv{}({\cal{O}},\pi)=\sum_{b\in {\cal{O}}(\pi)}\linv{}(b,\, \pi),\qquad 
\rinv{}({\cal{F}},\pi)=\sum_{b\in {\cal{F}}(\pi)}\rinv{}(b,\,\pi).\\
$$
\end{defi}

\begin{prop} Pour toute partition $\pi\in \P_n^k$, on a
\begin{eqnarray*}
\mak(\pi)&=&{\lmak}'(\pi)=\los(\pi)- \rinv{} ({\cal{F}},\pi)
+\linv({\cal{O}},\pi),\\
{\mak}'(\pi)&=&{\lmak}(\pi)=n(k-1)-\los(\pi)-\linv{} ({\cal{F}},\pi).
\end{eqnarray*}
\end{prop}
\pv  Soit
$\pi=B_1-\cdots -B_k \in P_{n}^{k}$, en utilisant les notations de
la d\'efinition~4,
 on peut r\'e\'ecrire
les statistiques $\lcs(\pi)$ et $\ros(\pi)$ de la fa\c con suivante~:
\begin{eqnarray*}
\lcs(\pi)&=&\sum_{i=1}^{k-1}\sum_{j>i}\left
(b_{j}-\rinv{}(g(B_{i}),B_{j})\right),\\
\ros(\pi)&=&\sum_{i=2}^{k}\sum_{j<i}\linv{}(p(B_{i}),B _{j}).
\end{eqnarray*}
Ainsi, la statistique $\mak(\pi)$ peut s'\'ecrire~:
\begin{eqnarray*}
\mak(\pi)&=&\sum_{i=1}^{k-1}\sum_{j>i}(b_{j}-
\rinv{}(g(B_{i}),B_{j}))+\sum_{i=2}^{k}\sum_{j<i}\linv{}(p(B_{i
}),B_{j}),\\ &=& \sum_{j=2}^k(j-1)b_j-\sum_{b\in
{\cal{F}}(\pi)}\rinv{}(b,\pi)+\sum_{b\in
{\cal{O}}(\pi)}\linv{}(b,\pi),\\ &=& \los(\pi)-\rinv{}
({\cal{F}},\pi) +\linv{} ({\cal{O}},\pi).
\end{eqnarray*}
D'autre part, on a~:
\begin{eqnarray*}
 \lcb(\pi)&=& \displaystyle \sum_{i=1}^{k-1}
\sum_{j>i}\rinv{}(g(B_{i}),B_{j})=\rinv{} ({\cal{F}},\pi),\\
\rob(\pi)&=&\sum_{i=2}^k\sum_{j<i}[b_j-\linv(p(B_i),B_j)]
=\sum_{j=1}^{k}(k-j)b_{j}-\linv{} ({\cal{O}},\pi).
\end{eqnarray*}
Et donc la statistique $\lmak$' peut s'\'ecrire~:
\begin{eqnarray*}
{\lmak}'(\pi)&=&\left[n(k-1)-\sum_{i=1}^{k}(k-i)b_{i}\right]-\rinv{}
({\cal{F}},\pi)+\linv{} ({\cal{O}},\pi)\\
&=&\left[ \sum_{i=1}^{k}(k-1)b_{i}-\sum_{i=1}^{k}(k-i)b_{i}\right] -\rinv{}
({\cal{F}},\pi)+\linv{} ({\cal{O}},\pi)\\
&=&\los(\pi)-\rinv{}({\cal{F}},\pi)+\linv ({\cal{O}},\pi).
\end{eqnarray*}
D'o\`u la premi\`ere identit\'e. La
seconde identit\'e peut \^etre v\'erifi\'ee de fa\c con analogue.\qed

%%%%%%%%%%%%%%%%%%%%%%%%%%%%%%%%%%%%%%%%%%%%%%%%%%%%%%%%%%%%%%%%
\section{Preuve du th\'eor\`eme~3}
%%%%%%%%%%%%%%%%%%%%%%%%%%%%%%%%%%%%%%%%%%%%%%%%%%%%%%%%%%%%%%%%
%%%
Dans tout ce qui suit on suppose que 
$\O$ est un sous-ensemble  fix\'e de $[n]$
\`a $k+1$ \'el\'ements avec $1\in \O$. 
Soit
$\P_{n}^{k+1}({\cal{O}})$ l'ensemble des partitions de $\P_{n}^{k+1}$
ayant pour l'ensemble des ouvrants ${\cal{O}}$. 

\begin{lem} La statistique $\los+\linv$
est constante sur $P_{n}^{k+1}({\cal O})$. Plus pr\'ecis\'ement, pour
toute partition $\pi\in P_{n}^{k+1}({\cal O})$, on a 
$$
\los(\pi)+\linv{} ({\cal{O}},\pi)=\sum_{x\in \O, x\neq 1}(n-x+1). 
$$
\end{lem}
\pv Soit $\pi_{0}=B_1-\cdots -B_{k+1}$ la partition de
$P_{n}^{k+1}({\cal O})$, telle que tout non-ouvrant $a$ soit le
plus \`a droite possible, c'est-\`a-dire, $a\in B_j$ tel que
$p(B_j)<a$ et $a<p(B_{j+1})$. 
Alors toutes les autres partitions $\pi$ de
$P_{n}^{k+1}({\cal O})$ s'obtiennent \`a partir de $\pi_{0}$ par
d\'eplacements successifs des non-ouvrants vers la gauche. Or
lorsque l'on d\'eplace une lettre vers la gauche (de $i$ blocs),
$\los(\pi_0)$ diminue de $i$, et $\linv{} ({\cal{O}},\pi_0)$
augmente de $i$. Ainsi on montre que $\los+\linv{}$ est constant
sur l'ensemble $P_{n}^{k+1}({\cal O})$. Il suffit donc calculer
 cette statistique pour $\pi_0$.  
Clairement
$\linv{}(\O,\pi_0)=0$. Si $(x_1,\ldots x_{k+1})$ est 
le r\'earrangement croissant des \'el\'ements de
$\O$, alors  pour tout $x_i\in \O$ le nombre d'\'el\'ements qui sont plus grand que $x_i$ et 
dans un bloc \`a droite de $B_i$  est $n-x_{i+1}+1$ pour $i=1,\ldots, k$.
D'o\`u 
le r\'esultat.
\qed

Pour tout $i\in
[k+1]$ et $\pi\in \P_{n}^{k+1}$ on pose
\begin{equation}\label{eq:statdef}
\stat{}_{i}(\pi)=k-\rinv{}({\cal{F}}, \pi)-\nrinv{}(g(B_{i}),\,
\pi).
\end{equation}
\begin{lem}
 Pour tout $i\in [k]$, il existe une bijection $\varphi_i$ sur
$\P_{n}^{k+1}({\cal{O}}) $ telle que~:
\begin{equation}\label{eq:stat}
\stat{}_{i}(\pi)=\stat{}_{i+1}(\varphi_i(\pi))-1.
\end{equation}
\end{lem}
\pv Soit $\pi = B_{1}-\cdots -B_{k+1}\in \P_{n}^{k+1}(\O)$. On
d\'efinit $\pi'=\varphi_i (\pi)=B_1'-\cdots -B_{k+1}'$ comme suit :
$B_j'=B_j$ pour tout $j\neq i,i+1$, $$
B_{i}'=\left\lbrace\begin{array}{ll} B_i\setminus \{a\in B_i\,|\,
a>g(B_{i+1})\}\cup \{g(B_{i+1})\} \qquad &\textrm{si}\quad
b_{i+1}>1;\\ B_i\setminus \{a\in B_i\,|\, a>g(B_{i+1})\} \qquad
&\textrm{si}\quad b_{i+1}=1;
   \end{array}\right.
$$ et $$ B_{i+1}'=\left\lbrace\begin{array}{ll} B_{i+1}\setminus
\{g(B_{i+1})\}\cup \{a\in B_i\,|\, a>g(B_{i+1})\} \qquad
&\textrm{si}\quad b_{i+1}>1;\\ B_{i+1}\cup \{a\in B_i\,|\,
a>g(B_{i+1})\} \qquad &\textrm{si}\quad b_{i+1}=1.
   \end{array}\right.
$$ 
On peut construire de mani\`ere analogue
l'application inverse $\varphi_i ^{-1}$. Donc $\varphi_i$ est une bijection.
Il reste \`a v\'erifier l'\'equation (\ref{eq:stat}). On distingue
trois cas suivants :
\begin{enumerate}
\item  Si $b_{i+1}=1$ et $g(B_i)<g(B_{i+1})$, alors $\pi'=\pi$.
\item  Si $b_{i+1}=1$ et $g(B_i)>g(B_{i+1})$ (et donc $b_i>1$), alors il est 
\'evident que
$\rinv{}({\cal{F}}, \pi')=\rinv{}({\cal{F}}, \pi)-1$, et
$\nrinv{}(g(B_{i+1}'), \pi')=\nrinv{}(g(B_{i}),\pi)$.
\item  Si $b_{i+1}>1$, on a
$\nrinv{}(g(B_{i+1}'),\pi')=b_{i+2}+\cdots +b_{k+1}
-\rinv(g(B_{i+1}'),\pi')$, et  $ \rinv{} ({\cal{F}},\pi')=\rinv{}
({\cal{F}},\pi)-\rinv(g(B_i),\pi)+b_{i+1}-1+\rinv(g(B_{i+1}'),\pi')$.
\end{enumerate}
Il est clair que l'\'egalit\'e (\ref{eq:stat}) a lieu dans les
deux premiers cas; pour le dernier cas, l'\'egalit\'e
(\ref{eq:stat}) d\'ecoule du fait que
$\rinv(g(B_i),\pi)+\nrinv(g(B_i),\pi)=b_{i+1}+\cdots + b_{k+1}$.
\qed

\begin{ex}
Prenons $i=3$, $k=3$, $n=9$, ${\O}=\{1,\, 2,\, 3, \, 5\}$ et
fixons $B_1=\{1,\, 4,\, 8\}$ et $B_2=\{2\}$.  Posons $\pi_{0} =
1\,4\,8-2-3-5\,6\,7\,9 $ et
 $\pi_{j}=\varphi_i(\pi_{j-1})$ pour $j\geq 1$, alors
l'application de $\varphi_i$ donne successivement~:
\begin{center}
\begin{tabular}{c|c|c}
                  partition                  &  $stat_{i+1}-1$  &
$stat_{i}$\\
 &&\\
\hline
  &&\\
       $\pi_{0} =  1\,4\,8-2-3-5\,6\,7\,9$   &      $-3$
&     $-6$\\
       $\pi_{1} =  1\,4\,8-2-3\,9-5\,6\,7$   &      $-6$
&     $-5$\\
       $\pi_{2} =  1\,4\,8-2-3\,7-5\,6\,9$   &      $-5$
&     $-5$\\
       $\pi_{3} =  1\,4\,8-2-3\,7\,9-5\,6$   &      $-5$
&     $-4$\\
       $\pi_{4} =  1\,4\,8-2-3\,6-5\,7\,9$   &      $-4$
&     $-5$\\
       $\pi_{5} =  1\,4\,8-2-3\,6\,9-5\,7$   &      $-5$
&     $-4$\\
       $\pi_{6} =  1\,4\,8-2-3\,6\,7-5\,9$   &      $-4$
&     $-4$\\
       $\pi_{7} =  1\,4\,8-2-3\,6\,7\,9-5$   &      $-4$
&     $-3$
\end{tabular}
\end{center}
On constate que $\pi_{j+8}=\pi_j$ pour $j\geq 0$.
\end{ex}

\begin{lem}
Pour tout $i \in \{1,\ldots,k\}$  on a
\begin{equation}\label{eq:lem}
\sum_{\pi \in
P_{n-1}^{k+1}}q^{\mak(\pi)+k-\nrinv{}(g(B_{i+1}), \pi)}=q^{i}\sum_{\pi \in 
P_{n-1}^{k+1}}q^{\mak(\pi)}.
\end{equation}
\end{lem}
\pv On montre d'abord par r\'ecurrence d\'ecroissante sur
$i$ l'identit\'e suivante~:
\begin{equation}\label{eq:ouv}
q^{i}\sum_{\pi \in P_{n-1}^{k+1}({\cal{O}})}q^{-\rinv{}
({\cal{F}},\pi)}=\sum_{\pi \in
P_{n-1}^{k+1}({\cal{O}})}q^{\stat{}_{i+1}(\pi)}.
\end{equation}
Pour $i=k$ le r\'esultat est vrai car
$\nrinv{}(g(B_{i+1}))=0$.
Supposons le r\'esultat vrai \`a
l'ordre $i$, alors 
$$ q^{i-1}\sum_{\pi \in
P_{n-1}^{k+1}({\cal{O}})}q^{-\rinv{} ({\cal{F}},\pi)}=\sum_{\pi \in
P_{n-1}^{k+1}({\cal{O}})}q^{\stat{}_{i+1}(\pi)-1}. $$ Or le
lemme~2 implique
$$
\sum_{\pi \in P_{n-1}^{k+1}({\cal{O}})}q^{\stat{}_{i+1}(\pi)-1}
=\sum_{\pi \in P_{n-1}^{k+1}({\cal{O}})}q^{\stat{}_{i}(\pi)}.
$$
Ce qui nous permet de conclure. Multipliant maintenant
les deux membres de  (\ref{eq:ouv}) 
par $q^{\sum_{x\in \O, x\neq 1}(n-x+1)}$, on 
d\'eduit du lemme~4~: 
$$
 q^{i}\sum_{\pi \in P_{n-1}^{k+1}({\cal
O})}q^{\linv{}({\cal{O}},\pi)+\los(\pi)-\rinv{}(\cal{F},\pi)}=\sum_{ \pi \in
P_{n-1}^{k+1}({\cal
O})}q^{\stat{}_{i+1}(\pi)+\linv{}({\cal{O}},\pi)+ \los(\pi)}. $$
Compte tenu de (\ref{eq:statdef}) et la proposition~1, on obtient (13) en
sommant sur tous les ouvrants  $\O$ possibles.
\qed

On est maintenant en mesure de d\'emontrer que la fonction g\'en\'eratrice 
de $\mak_l$ sur $\P_n^k$ v\'erifie la relation de 
r\'ecurrence~(\ref{eq:stirling})
On d\'emontre d'abord  ce r\'esultat 
pour $\mak$ sur ${\P}_n^k$.
 Il est clair que c'est vrai pour
$n=1$. Etant donn\'ee une partition $\pi\in {\P}_n^{k+1}$, on note
$\pi'$ la partition obtenue en supprimant $n$.
 On
distingue alors deux cas selon que $n$ est un singleton ou non.
\begin{enumerate}
\item $n$ est un singleton.
Alors $\pi'\in P_{n-1}^{k}$ et on v\'erifie sans peine que
\begin{eqnarray*}
\los(\pi)&=&\los(\pi')+k,\\ \linv{} ({\cal{O}},\pi)&=&\linv{}
({\cal{O}},\pi'),\\ \rinv{} ({\cal{F}},\pi)&=&\rinv{}
({\cal{F}},\pi').
\end{eqnarray*}
Ainsi $\mak(\pi)=\mak(\pi')+k$. D'o\`u on d\'eduit
 la fonction
g\'en\'eratrice correspondante~: $$ q^{k}\sum_{\pi'\in
{\P}_{n-1}^{k}}q^{\mak(\pi')}. $$
\item $n$ n'est pas un singleton. Alors
$\pi'\in P_{n-1}^{k+1}$. Supposons que $n$ soit dans le
$i^\textrm{\`eme}$ bloc de $\pi$. Alors
\begin{eqnarray*}
\los(\pi)&=&\los(\pi')+i-1,\\ \linv{} ({\cal{O}},\pi)&=&\linv{}
({\cal{O}},\pi')+k+1-i,\\ \rinv{} ({\cal{F}},\pi)&=&\rinv{}
({\cal{F}},\pi')-\rinv{}(g(B_{i}),\pi)+b_{i+1}+\cdots +b_{k+1}.
\end{eqnarray*}
Ainsi, en vertu de la proposition~1, on a $$\mak(\pi)=\mak(\pi')+k
-\nrinv{}(g(B_{i}),\pi). $$ On en d\'eduit donc, d'apr\`es le
lemme~4, la fonction g\'en\'eratrice correspondante~: 
$$\sum_{i=0}^{k} \sum_{\pi
\in
P_{n-1}^{k+1}}q^{\mak(\pi)+k-\nrinv{}(g(B_{i+1}),\pi)}=[k]_q\sum_{\pi \in
P_{n-1}^{k+1}}q^{\mak(\pi)}. $$
\end{enumerate}

En r\'ecapitulant les deux cas pr\'ec\'edants, on a~:
\begin{eqnarray*}
\sum_{\pi \in P_{n}^{k+1}}q^{\mak(\pi)} &=&q^{k} \sum_{\pi \in
P_{n-1}^{k}}q^{\mak(\pi)}+[k]_q\sum_{\pi \in
P_{n-1}^{k+1}}q^{\mak(\pi)}.
\end{eqnarray*}
Ce qui est exactement la relation de r\'ecurrence~(\ref{eq:stirling})
pour les
$q$-nombres de Stirling $S_q(n,k)$. Donc la fonction
g\'en\'eratrice de $\mak$ sur $\P_n^k$ est $S_q(n,k)$.

Enfin, pour tout $l\in [k]$, l'\'equation ~(\ref{eq:lem}) \'equivaut
\`a :
$$
\sum_{\pi \in P_{n}^{k}}q^{\mak(\pi)}=\sum_{\pi \in
P_{n}^{k}}q^{k-l+\mak(\pi)-\nrinv{}(g(B_l), \pi)}.
$$
Ce qui montre que $\mak$ et $\mak_l$ sont \'equidistribu\'ees sur  $\P_n^k$.

%%%%%%%%%%%%%%%%%%%%%%%%%%%%%%%%%%%%%%%%%%%%%%%%%%%%%%%%%%%%
\section{Remarques sur les partitions ordonn\'ees}
%%%%%%%%%%%%%%%%%%%%%%%%%%%%%%%%%%%%%%%%%%%%%%%%%%%%%%%%%%
Une 
{\sl $k$-partition ordonn\'ee} de $[n]$ est 
une suite  $(B_1,\, B_2,\,
\ldots, B_k)$ de $k$ sous-ensembles de $[n]$ telle que 
$\pi=B_{\sigma(1)}-B_{\sigma(2)}-\cdots -B_{\sigma(k)}$ soit une partition de $\P_n^k$ pour une 
permutation $\sigma$ de $[k]$.
Notons $\OP_n^k$ l'ensemble des $k$-partitions ordonn\'ees de  $[n]$.
Il est \'evident que le cardinal de 
$\OP_n^k$ est $k!S_1(n,k)$. Il s'agit de trouver des statistiques 
{\sl Euler-mahoniennes}
sur $\OP_n^k$, i.e., leurs fonctions g\'en\'eratrices sur $\OP_n^k$
 sont \'egales \`a $[k]_q!S_q(n,k)$. Certaines de ces statistiques 
peuvent \^etre obtenues \`a partir d'un r\'esultat de Wachs~\cite{Wa}.
Steingr\'{\i}msson\cite{St} en a  propos\'e
d'autres.

\begin{defi}
Soit $\pi=B_1-\cdots-B_k \in \OP_n^k$, on d\'efinit un ordre
partiel sur les blocs comme suit~: $B_i>B_j$ si toutes
les lettres de $B_i$ sont plus grandes que celles de $B_j$.
 On dit que $i$ est un indice de descente si 
$B_{i}>B_{i+1}$.
On d\'efinit alors $\bmaj(\pi)$ comme la somme de tous les indices de
descentes de $\pi$; et $\binv(\pi)$ comme le nombre de couples $(i,j)$ 
tel que $i<j$ et $B_{i}>B_{j}$.
\end{defi}

Steingr\'{\i}msson\cite[Conj.~13]{St} a conjectur\'e que si l'on ajoute
\`a $\bmaj$ ou $\binv$, l'une des statistiques de la d\'efinition 1,
on obtient une statistique  Euler-mahonienne, c'est \`a dire~:
\begin{conj}[Steingr\'{\i}msson] Les
statistiques suivantes sont  Euler-mahoniennes~:\\
$$
\begin{tabular}{llll}
$\mak+\bmaj$, & $\qquad \mak{}'+\bmaj$, & $\qquad \lmak{}'+\bmaj$, & $\qquad \lmak+\bmaj$,\\
$\mak+\binv$, & $\qquad \mak{}'+\binv$, & $\qquad \lmak{}'+\binv$, & $\qquad \lmak+\binv$.
\end{tabular}
$$
\end{conj}

On remarque que la d\'emonstration  de la proposition~1
 s'\'etend \emph{mutatis mutandis} au cas des partitions de $OP_{n}^{k}$, 
on peut alors r\'eduire cette conjecture de moiti\'e, i.e., 
dans la conjecture ci-dessus, il n'y a que quatre statistiques distinctes. Plus pr\'ecis\'ement 
on a le r\'esultat suivant~:
\begin{prop} On a les \'egalit\'es suivantes~:
\begin{eqnarray*}
\mak+\bmaj&=&\lmak{}'+\bmaj,\qquad \mak{}'+\bmaj=\lmak+\bmaj,\\
\mak+\binv&=&\lmak{}'+\binv, \qquad\hspace{8pt} \mak{}'+\binv=\lmak+\binv.
\end{eqnarray*}
\end{prop}

%%%%%%%%%%%%%%%%%%%%%%%%%%%%%%%%%%%%%%%%%%%%%%%%%%%%%%%%%%%%%%%%
\medskip 
\parindent=0pt {\bf Remerciements}
Les auteurs tiennent \`a remercier les deux arbitres anonymes pour leurs
lectures attentives sur une version ant\'erieure, permettant
d'am\'eliorer la r\'edaction de cet article. 
%%%%%%%%%%%%%%%%%%%%%%%%%%%%%%%%%%%%%%%%%%%%%%%%%%%%%%%%%%%%%%%%
%%%%%


\begin{thebibliography}{9}
\bibitem{CA}  Carlitz (L.),
\emph{On Abelian fields}, Trans. Amer. Math. Soc.,
{\bf 35} (1933), 122-136.
%%%%%%%%%%%%%%%%%%%%%%%%%%%%%%%%%%%%%%%%%%%%%%%%%%%%%%%%%%%%%%%%
%%%%%%%
\bibitem{DL} De M\'edicis (A.) and Leroux (P.), \emph{A unified
combinatorial approach for $q$-(and $p,q$-) Stirling numbers}, J.
Statistical Planning and Inference, {\bf 34} (1993), 89-105.

%%%%%%%%%%%%%%%%%%%%%%%%%%%%%%%%%%%%%%%%%%%%%%%%%%%%%%%%%%%%%%%%
%%%%%%%
\bibitem{EH} Ehrenborg (R.) and Readdy (M.),
\emph{Juggling and applications to q-analogues}, Discrete Mathematics,
Special Issue on Algebraic Combinatorics, {\bf 157} (1996), 107-125.
%%%%%%%%%%%%%%%%%%%%%%%%%%%%%%%%%%%%%%%%%%%%%%%%%%%%%%%%%%%%%%%%
\bibitem{Fl}  Flajolet (Ph.),
\emph{Combinatorial aspects of continued fractions}, Disc. Math.,
{\bf 41}(1982), 125-161.
%%%%%%%%%%%%%%%%%%%%%%%%%%%%%%%%%%%%%%%%%%%%%%%%%%%%%%%%%%%%%%%%
%%%%%%%
\bibitem{FZ}  Foata (D.) and  Zeilberger (D.),
\emph{Denert's permutation statistic is
    indeed Euler-Mahonian}, Studies in Appl. Math., {\bf 83} (1990),
31-59.
%%%%%%%%%%%%%%%%%%%%%%%%%%%%%%%%%%%%%%%%%%%%%%%%%%%%%%%%%%%%%%%%
%%%%%%
\bibitem{GR} Garsia (A.) et Remmel (J.B.),
\emph{$Q$-counting rook configurations and a formula of
Frobenius}, J. Combin. Theory Ser. A, {\bf 41} (1986), 246--275.
%%%%%%%%%%%%%%%%%%%%%%%%%%%%%%%%%%%%%%%%%%%%%%%%%%%%%%%%%%%%%%%%
\bibitem{GO} Gould (H.W.),
\emph{The $q$-Stirling numbers of the first and second kinds},
Duke Math. J., {\bf 28} (1961), 281-289.
%%%%%%%%%%%%%%%%%%%%%%%%%%%%%%%%%%%%%%%%%%%%%%%%%%%%%%%%%%%
\bibitem{LE} Leroux (P.),
\emph{Reduced Matrices and q-log-concavity Properties of q-stirling Numbers},
J. Combin. Theory Ser. A, {\bf 54} (1990), 64-84.
%%%%%%%%%%%%%%%%%%%%%%%%%%%%%%%%%%%%%%%%%%%%%%%%%%%%%%%%%%%
 \bibitem{Mi}  Milne (S.), \emph{Restricted growth functions, rank row
matching of partition lattices, and q-stirling numbers},
Adv. Math., {\bf 43} (1982), 173-196.
%%%%%%%%%%%%%%%%%%%%%%%%%%%%%%%%%%%%%%%%%%%%%%%%%%%%%%%%%%%
\bibitem{Sa}  Sagan (B.), \emph{A maj statistics for set
partitions}, European J. Combin., {\bf 12}~(1991), 69-79.
%%%%%%%%%%%%%%%%%%%%%%%%%%%%%%%%%%%%%%%%%%%%%%%%%%%%%%%%%%%
\bibitem{St}  Steingr\'{\i}msson (E.), \emph{Statistics on
ordered partitions of  sets},  preprint, 1999.

%%%%%%%%%%%%%%%%%%%%%%%%%%%%%%%%%%%%%%%%%%%%%%%%%%%%%%%%%%%
\bibitem{Wa} Wachs ( M.), \emph{$\sigma$-Restricted Growth Functions and 
p,q-stirling
numbers}, J. Combin. Theory Ser. A, {\bf 68} (1994), 470-480.
%%%%%%%%%%%%%%%%%%%%%%%%%%%%%%%%%%%%%%%%%%%%%%%%%%%%%%%%%%%
\bibitem{WW} Wachs ( M.) and  White (D.), \emph{p,q-stirling
numbers and set partition statistics}, J. Combin. Theory Ser. A,
{\bf 56} (1991), 27-46.
%%%%%%%%%%%%%%%%%%%%%%%%%%%%%%%%%%%%%%%%%%%%%%%%%%%%%%%%%%%%%%%
\bibitem{Wh} White (D.), \emph{Interpolating Set Partition Statistics},
J. Combin. Theory Ser.~A, {\bf 68} (1994), 262-295.
%%%%%%%%%%%%%%%%%%%%%%%%%%%%%%%%%%%%%%%%%%%%%%%%%%%%%%%%%%%%%%%
\end{thebibliography}
\end{document}